\documentclass{article}

\usepackage{latexsym,amsmath,amssymb,amsthm,bm,booktabs,authblk,multicol}
\usepackage{setspace}
\usepackage{epsf}

\newcommand{\bbR}{\mathbb{R}}

\newcommand{\bfx}{\boldsymbol{x}}

\newcommand{\bfc}{\boldsymbol{c}}
\newcommand{\bfd}{\boldsymbol{d}}

\newcommand{\bfb}{\boldsymbol{b}}

\newcommand{\bfg}{\boldsymbol{g}}
\newcommand{\bfzero}{\boldsymbol{0}}

\newcommand{\minimize}{\text{minimize}}
\newcommand{\st}{\text{subject to}}

\newtheorem{thm}{Theorem}
\newtheorem{lem}{Lemma}
\newtheorem{cor}{Corollary}

\qedhere

\begin{document}
%\renewcommand\arraystretch{1.1}
%\renewcommand\baselinestretch{1.1}

%%%

\title{A primal-simplex based Tardos' algorithm}

\date{\empty}

\author[a]{Shinji Mizuno}
\author[a]{Noriyoshi Sukegawa} %\thanks{e-mail: sukegawa.n.aa@m.titech.ac.jp}
\author[b]{Antoine Deza}

\affil[a]{{\small Graduate School of Decision Science and Technology, Tokyo Institute of Technology, 
2-12-1-W9-58, Oo-Okayama, Meguro-ku, Tokyo, 152-8552, Japan. }}
\affil[b]{{\small Advanced Optimization Laboratory, Department of Computing and Software, McMaster University, Hamilton, Ontario, Canada.}}

\maketitle
\begin{abstract}
\noindent
In the mid-eighties Tardos proposed a strongly polynomial algorithm for solving linear programming problems
for which the size of the coefficient matrix is polynomially bounded by the dimension. 
Combining Orlin's primal-based modification and Mizuno's use of the simplex method, 
we introduce a modification of Tardos' algorithm considering only the primal problem 
and using simplex method to solve the auxiliary problems. The proposed algorithm is strongly 
polynomial if  the coefficient matrix is totally unimodular and the auxiliary problems are non-degenerate. 
\end{abstract}

\noindent
%\begin{center}
{\bf Keyword}:~Tardos' algorithm, simplex method, strongly polynomial algorithm, total  unimodularity
%\end{center}

%%%
\section{Introduction} \label{sec:Intro}
%%%

In the mid-eighties Tardos~\cite{Ta85,Ta86} proposed a strongly polynomial algorithm 
for solving linear programming problems $\min \{ \bfc^\top \bfx \,|\, A\bfx=\bfb, \, \bfx \geq \bfzero \}$
for which the size of the coefficients of $A$ are polynomially bounded by the dimension. 
Such instances include minimum cost flow, bipartite matching, multicommodity flow, and vertex packing in chordal graphs.
The basic strategy of Tardos' algorithm is to identify the coordinates equal to zero at optimality. 
The algorithm involves solving several auxiliary dual problems by the ellipsoid or interior-point methods. 
By successively identifying such vanishing coordinates, the problem is made smaller and an optimal solution 
is obtained inductively.  Orlin~\cite{Or86} proposed a modification of Tardos' algorithm considering only the 
primal problem; that is, identifying the coordinates strictly positive at optimality. 
He observed that the right-hand side coefficients of the auxiliary problems might be impractically large.

\par 
In 2014, Mizuno~\cite{Mi14} modified Tardos' algorithm by using a dual simplex method to solve the auxiliary problems. 
He observed that this approach is strongly polynomial if $A$ is totally unimodular and  the auxiliary problems are non-degenerate;
that is, the basic variables are strictly positive for every basic feasible solution. 
The strong polynomiality is a consequence of Kitahara and Mizuno~\cite{KiMi12, KiMi13} results
which extend in part Ye's  result~\cite{Ye11} for Markov decision problems  and bounds  the 
number of distinct basic feasible solutions generated by the simplex method. 

\par
Combining Orlin's and Mizuno's approaches, 
we introduce a modification of the algorithm proposed by Mizuno considering only the primal problem. 
The proposed algorithm is strongly polynomial if $A$ is totally unimodular and the auxiliary problems are non-degenerate. 
As it involves only the primal and does not suffer from impractically large right-hand side coefficients, 
the  proposed algorithm improves the implementability of the approach. 
While the proposed algorithm and the complexity analysis is focusing on the case where $A$ is totally unimodular, the algorithm could
be enhanced  to handle general matrices. 
The enhanced algorithm would be strongly polynomiality if the absolute 
value of any subdeterminant of $A$ is polynomially bounded by the dimension.

%%%
\section{A primal-simplex based Tardos' algorithm} \label{sec:Alg}
%%%

\subsection{Formulation and main result}

Consider the following formulation:
\begin{equation}\label{P}
\begin{array}{lllllllllllll}
\minimize	&\bfc^\top \bfx						\\	
\st 			&A \bfx = \bfb, \, \bfx \geq \bfzero	 
\end{array}
\end{equation}
where $A \in \bbR^{m \times n}$, $\bfb \in \bbR^{m}$, and $\bfc \in \bbR^{n}$ are given. 
% In what follows, assume that $A$ is totally unimodular. 
% We denote by $\Delta$ the largest absolute value of the sub determinant of $A$; that is,
%  We recall that a matrix $A$ is totally unimodular if all the sub-determinants are equal to either -1, 0 or 1.
% In what follows, we focus on the case where $A$ is totally unimodular. 
% Therefore, if we let 
% \begin{equation*}
% \Delta = \max \{\, {|{\det D}|} \, | \, D \text{ is a square submatrix of } A \,\}, 
% \end{equation*}
% we have $\Delta =1$.
The optimal solution of $(\ref{P})$, if any, is assumed without loss of generality to be unique.
Otherwise $\bfc$ could be perturbed by $(\epsilon, \epsilon^2,\ldots, \epsilon^n)$ for a sufficiently small $\epsilon>0$.
Alternatively, the simplex method can be performed using a lexicographical order if a tie occurs when choosing an 
entering variable by Dantzig's rule. Let $K^* \subseteq N=\{1, 2, \ldots, n \}$ be the optimal basis of $(\ref{P})$. 
The proposed algorithm inductively builds a subset $\bar{K} \subseteq K^*$ through solving an auxiliary problem. 
If $\bar{K} = K^*$ we obtained the optimal solution.  Otherwise, we obtain a smaller yet equivalent problem 
by deleting the variables corresponding to $\bar{K}$.  Thus, the optimal solution is obtained inductively. 
For clarity of the exposition of the algorithm and of the proof of Theorem~\ref{thm1}, we assume in the remainder 
of the paper that $A$ is totally unimodular; that is,  all its subdeterminants are equal to either $-1, 0$ or $1$.

\begin{thm}\label{thm1}
The primal-simplex based Tardos' algorithm is strongly polynomial if $A$ is totally unimodular and all the auxiliary problems are non-degenerate;
that is, all the basic variables are strictly positive for every basic feasible solution. 
\begin{proof}
See Section~\ref{sec:Proof}. 
\end{proof}
\end{thm}
\noindent
%The detail of our algorithm is as follows. 

\subsection{A primal-simplex based Tardos' algorithm}

%The algorithm is stated as follows: 
\begin{description}

\item[Step 0 (initialization):] $\:$\\\\
Let $\bar{K} : = \emptyset$ and its complement $K : = N$. 
%Set $\bfc' = \bfc$, $\bfb' = \bfb$, and $A' = A$. 
% Go to Step 1. 

\item[Step 1 (reduction):] $\:$\\\\
If $\bar{K} \ne \emptyset$, remove the variables corresponding to $\bar{K}$ in the following way.\\
Let $G \in \bbR^{m \times m}$ be a nonsingular submatrix of $A$ such that 
its first $|\bar{K}|$ columns form $A_{\bar{K}}$ and $H=G^{-1}$.  Let $H_1$ consists of the 
first $|\bar{K}|$ rows of $H$, $H_2$ denote the remainder, and 
consider the following {\em reduced} problem: 
\begin{equation}\label{smallP}
	\begin{array}{ll}
		\minimize	&\bfc'^\top \bfx'		\\
		\st 			&A' \bfx' = \bfb', \, \bfx' \geq \bfzero,  	 
	\end{array}
\end{equation}
where $A' = H_2 A_K$, $\bfb' = H_2 \bfb$, $\bfc' = \bfc_K - (H_1 A_K)^\top \bfc_{\bar{K}}$, and $\bfx' = \bfx_K$. \\%\\\\
If $\bar{K} = \emptyset$, set $A' : = A$, $\bfb' : = \bfb$, and $\bfc' : = \bfc$. \\
 {\sc Go to Step 2}.
%When $\bar{K} \ne \emptyset$, we remove the variables $x_i$ $(i \in \bar{K})$ from the problem (\ref{P}) 
%by the Gaussian elimination, that is, express $\bfx_{\bar{K}}=\bar{\bfb} - \bar{A}_{K} \bfx_{K}$  
%by using $|\bar{K}|$ equalities in $A \bfx = \bfb$ and  
%substitute $\bar{\bfb} - \bar{A}_{K} \bfx_{K}$ for $\bfx_{\bar{K}}$ in the other equalities and 
%the objective function. 
%Then we obtain the following ``reduced'' problem: 
%\begin{equation}\label{smallP}
%\begin{array}{lllllllllllll}
%\minimize	&\bfc'^\top \bfx'		\\
%\st 			&A' \bfx' = \bfb', \, \bfx' \geq \bfzero,  	 
%\end{array}
%\end{equation}
%where $\bfx' = \bfx_{K}$  
%and $A'$ is an $m' \times n'$ matrix for $m'=m-|\bar{K}|$ and $n'=n-|\bar{K}|$. 
%When $\bar{K} = \emptyset$, we have the problem $(\ref{smallP})$ for $\bfc' = \bfc$, $\bfb' = \bfb$, and $A' = A$.  
%Go to Step 2. 

\item[Step 2 (scaling and rounding):] $\:$\\\\
Let $m'=m-|\bar{K}|$ and $n'=n-|\bar{K}|$. % denote the number of rows and columns of $A'$. 
For a basis $L \subseteq K$ of $A'$ and $\bar{L} = K \setminus L$, rewrite $(\ref{smallP})$ as: 
\begin{equation}\label{BP}
\begin{array}{lllllllllllll}
\minimize	&\bfc'^\top \bfx'												\\	
\st 			&\bfx'_L + (A_{L}')^{-1} A'_{\bar{L}} \bfx'_{\bar{L}} = (A_{L}')^{-1} \bfb', 
\, \bfx' \geq \bfzero.	 
\end{array}
\end{equation}
If $(A_{L}')^{-1} \bfb' = \bfzero$, stop. 
Otherwise, consider the following {\em scaled}  problem: 
\begin{equation}\label{scaleP}
\begin{array}{lllllllllllll}
\minimize	&{\bfc'}^\top \bfx'											\\	
\st 			&\bfx'_L + (A_{L}')^{-1} A'_{\bar{L}} \bfx'_{\bar{L}} = (A_{L}')^{-1} \bfb'/k, 
\, \bfx' \geq \bfzero, 
\end{array}
\end{equation}
where $k={\| {A'}^\top ({A'}{A'}^\top)^{-1}\bfb' \|_2}/{(m' + (n')^2)}$. Then, consider the following {\em rounded} problem: 
\begin{equation}\label{relaxP}
\begin{array}{lllllllllllll}
\minimize	&{\bfc'}^\top \bfx'											\\	
\st 			&\bfx'_L + (A_{L}')^{-1} A'_{\bar{L}} \bfx'_{\bar{L}} = \lceil (A_{L}')^{-1} \bfb' / k \rceil, 
\, \bfx' \geq \bfzero.	 
\end{array}
\end{equation}
If $(\ref{relaxP})$ is infeasible, stop. 
Otherwise, solve $(\ref{relaxP})$ using the simplex method with Dantzig's rule. 
If $(\ref{relaxP})$ is unbounded, stop.  Otherwise, let $\bfx''$ denote the optimal solution 
and $L''$ the optimal basis. If $\bar{K} \cup L''$ is an optimal basis of the original problem (\ref{P}), stop. 
Otherwise, {\sc go to Step 3}. 

\item[Step 3 (iteration):] $\:$\\\\
Set  $\bar{K} : = \bar{K} \cup J$ and $K : = K \setminus J$ where  $J = \{\, i \,|\, x''_i \ge n', \, i \in K \,\}$.\\
If $|K|=n-m$, stop.  Otherwise, {\sc go to Step 1}. 
\end{description}

\subsection{Annotations of  the proposed algorithm}

If $\bar{K} \subseteq K^*$ and the optimal solution of $(\ref{P})$ is unique, 
we can remove the non-negativity constraints for  $x_i$ for $i \in \bar{K}$. 
In Step 1, the reduced problem $(\ref{smallP})$ is obtained by
 expressing $\bfx_{\bar{K}}$ as $H_1 \bfb - H_1 A_K \bfx_K$
and substituting $H_1 \bfb - H_1 A_K \bfx_K$ for $\bfx_{\bar{K}}$ in the objective function.
Therefore, the optimal solution for $(\ref{smallP})$ yields  the optimal solution for $(\ref{P})$ 
via $\bfx_{\bar{K}} = H_1\bfb - H_1 A_K \bfx_{K}$. 
%We observe that the optimal solution of the reduced problem $(\ref{smallP})$ gives us that of $(\ref{P})$ 
%via $\bfx_{\bar{K}} = \bar{\bfb} - \bar{A}_{K} \bfx_{K}$ since $\bar{K} \subseteq K^*$. 
The constant term in the objective function is removed for simplicity. 
Note that the matrices $A'$ and $[ I,\, (A_{L}')^{-1} A'_{\bar{L}} ]$ involved in
$(\ref{smallP})$,  $(\ref{BP})$, $(\ref{scaleP})$, and $(\ref{relaxP})$ are totally unimodular
 if $A$ is totally unimodular,  see Theorem 19.5 in Schrijver~\cite{Sc86}.\\

\par
%If we have $(A_{L}')^{-1} \bfb' = \bfzero$ in Step 2, 
%$\bfx'=\bfzero$ is the optimal solution of the reduced problem $(\ref{smallP})$ 
%unless the problem is unbounded. Note that our algorithm never stops here 
%under nondegenerate assumption.  
In step 2, the scaling factor $k$ is strictly positive if $(A_{L}')^{-1} \bfb' \neq \bfzero$ and, see Lemma~\ref{lem:gamma},
$\| \lceil (A_{L}')^{-1} \bfb' / k \rceil \|_\infty$ is polynomially bounded above in $m'$ and $n'$,
which is a key fact for showing the strong polynomiality.  
Although the proposed algorithm builds the simplex tableau associated to $(\ref{BP})$ and the reduced problem $(\ref{smallP})$ from
scratch at each iteration, it is essentially for clarity of the exposition and can be ignored. In particular, one can observe that 
$L'' \setminus J$ can be used as the basis $L$ for (\ref{BP}) at the next iteration, thus enabling a warm start. 
%(This bound becomes smaller as we proceed the algorithm 
%because $m'=m-|\bar{K}|$, $n'=n-|\bar{K}|$, and $|\bar{K}|$ is increased (at least) by one for each iteration. 
%(S) I can not observe this statement, so some explanation is needed. 
%(N) Now I come to think that this part is tiny thing. )
%Note that Orlin~\cite{Or86} uses huge (exponentially large) values to this part, which makes his algorithm impractical. 
By performing Phase one of the two-phase simplex method for the rounded problem $(\ref{relaxP})$, 
we can check the feasibility of $(\ref{relaxP})$ and compute an initial basic feasible solution, unless it is infeasible.\\  
%The feasibility of the rounded problem $(\ref{relaxP})$ is checked via solving several auxiliary LP problems,
%and yields  an initial basis, see Section~\ref{sec:Complexity}. 
%Note that the infeasibility of the rounded problem $(\ref{relaxP})$ implies the infeasibility of the original problem $(\ref{P})$ 
%as well as the scaled problem $(\ref{scaleP})$. 
%Also note that the original problem $(\ref{P})$ is either unbounded or infeasible when the rounded problem $(\ref{relaxP})$ is unbounded. 

\par
In Step 3,  $J\neq \emptyset$ by Lemma~\ref{lem:Scale}; that is,  the size of $K$ is strictly decreasing. Thus, the proposed 
algorithm terminates after at most $m$ iterations. %a finite number of iterations. 
If $(\ref{P})$ has an optimal solution, $\bar{K} \subseteq K^*$ by Corollary~\ref{lem:Reduction}.\\

\par The stopping conditions of the proposed  algorithm are: 
\begin{itemize}
\setlength{\itemsep}{0em}
\item[$\circ$] if $(A_{L}')^{-1} \bfb' = \bfzero$,  
the simplex tableau associated to $(\ref{BP})$ yields either the optimality of $\bfx'=\bfzero$ or the unboundedness of the reduced problem $(\ref{smallP})$. 
%We can check the boundedness by solving the problem  $(\ref{smallP})$ directly by the simplex method.  
%In the former case, we obtain the optimal solution of the problem $(\ref{P})$ 
%via $\bfx_{\bar{K}} = \bar{\bfb} - \bar{A}_{K} \bfx_{K}$. 
%In the latter case, $(\ref{P})$ is also unbounded. 
\item[$\circ$]
since the rounded problem $(\ref{relaxP})$ is a relaxation of the scaled problem $(\ref{scaleP})$,
\begin{itemize}
\setlength{\itemsep}{0em}
\item
the scaled problem $(\ref{scaleP})$ and the original problem $(\ref{P})$ are both infeasible if $(\ref{relaxP})$ is infeasible
\item
 the scaled problem  $(\ref{scaleP})$ is unbounded or infeasible if $(\ref{relaxP})$ is unbounded. In both
cases, the original problem $(\ref{P})$ has no optimal solution. 
\end{itemize}
\item[$\circ$]
if $|K|=n-m$ in Step 3, the problem $(\ref{P})$ is infeasible as otherwise the algorithm finds an optimal basis in Step 2. 
\end{itemize}

\section{Proof of Theorem~\ref{thm1}}\label{sec:Proof}

Lemma~\ref{lem:Scale} states that  the set $J = \{ i \,|\, x''_i \ge n', \, i \in K \}$ used  in Step 3
 is never empty and thus,  the proposed algorithm solves the rounded problem $(\ref{relaxP})$ at most $m$ times. 
\begin{lem}\label{lem:Scale}
 $J\neq \emptyset$  as any solution $\bfx''$ of the rounded problem $(\ref{relaxP})$ 
satisfies $\| \bfx'' \|_{\infty} \ge n'$.
\begin{proof}
Let  $\bfx''$ be any solution of the rounded problem $(\ref{relaxP})$. Then 
\[
  A' \bfx'' = A'_L \lceil (A_{L}')^{-1} \bfb' / k \rceil. 
\]
Since, for  any $\bfg$, the minimal $l_2$-norm point satisfying  $A' \bfx' = \bfg$ is $A'^T(A' A'^T)^{-1} \bfg$, we have
\[ \begin{array}{lll}
  \| \bfx'' \|_2 & \ge & \|A'^T(A' A'^T)^{-1} A'_L \lceil (A_{L}')^{-1} \bfb' / k \rceil \|_2 \\
  & \ge &  \|A'^T(A' A'^T)^{-1} \bfb' / k \|_2 -
    \|A'^T(A' A'^T)^{-1} A'_L \bfd \|_2 \\
  & = & (m' + (n')^2 ) -   \|A'^T(A' A'^T)^{-1} A' 
   \left( \begin{array}{l} \bfd \\ \bfzero_{\bar L} \end{array} \right) \|_2 \\
  & \ge & (n')^2 + m'-   \|\bfd \|_2, 
 \end{array}
\]
where $k={\| {A'}^\top ({A'}{A'}^\top)^{-1}\bfb' \|_2}/{(m' + (n')^2)}$ and
$\bfd = (A_{L}')^{-1} \bfb' / k - \lceil (A_{L}')^{-1} \bfb' / k \rceil$. 
Since $\| \bfd \|_{\infty} < 1$, we obtain that 
\[
    \| \bfx'' \|_{\infty} \ge  \| \bfx'' \|_2/n'  >  ((n')^2 + m' - m')/n' = n'.  
\] 
%
% 
%denote the feasible region of the rounded problem $(\ref{relaxP})$ without the non-negativity constraints. 
%Let $E=[ I,\, (A_{L}')^{-1} A'_{\bar{L}} ]$ and $X=E^\top (EE^\top)^{-1}$. 
%The minimal $l_2$-norm point of the system $E \bfx' = \lceil (A_{L}')^{-1} \bfb' / k \rceil$ 
%is determined  as $X \lceil (A_{L}')^{-1} \bfb' / k \rceil$ and, by the triangle inequality, 
%\begin{equation*}
%\| X \lceil (A_{L}')^{-1} \bfb' / k \rceil \|_2 \ge \|X(A_{L}')^{-1} \bfb' \|_2 / k - \|  X \bfd \|_2
%\end{equation*}
%where $\bfd = (A_{L}')^{-1} \bfb' / k - \lceil (A_{L}')^{-1} \bfb' / k \rceil$. 
%Note that $X(A_{L}')^{-1} \bfb'$ is the minimal $l_2$-norm point of the system $A' \bfx' = \bfb'$. 
%Then, $\|X(A_{L}')^{-1} \bfb' \|_2 / k = (m')^2 + (n')^2$ by the definition of $k$.
%
%\par
%Since $X \bfd$ is the minimal $l_2$-norm point of the system $E \bfx' = \bfd$, we have
%$\| X \bfd \|_2 \le \| \bfy' \|_2 \le m' \|\bfy' \|_\infty $ for any basic solution $\bfy'$ of the system $E \bfx' = \bfd$. 
%Since $E$ is totally unimodular and $\| \bfd \|_\infty < 1$, 
%we have $\| \bfy' \|_\infty \le m'$  which implies $\| X \lceil (A_{L}')^{-1} \bfb' / k \rceil \|_2 \ge (n')^2$. 
%This completes the proof as  $\bfx''$ satisfies $E \bfx'' = \lceil (A_{L}')^{-1} \bfb' / k \rceil$, 
%and $\| \bfx'' \|_\infty \ge \| \bfx'' \|_2 / n'$. 
\end{proof}
\end{lem}
Corollary~\ref{lem:Reduction} is a direct consequence of Theorem~\ref{lem:Sc}
and shows that  $\bar{K} \subseteq K^*$. 

\begin{thm}[Theorem 10.5 in Schrijver~\cite{Sc86}]\label{lem:Sc}
Let $A$ be an $m \times n$-matrix, and let $\Delta^*$ be such that for each nonsingular submatrix 
$B$ of $A$ all entries of $B^{-1}$ are at most $\Delta^*$ in absolute value. 
Let $\bfc$ be a column $n$-vector, and let $\bfb''$ and $\bfb^*$ be column $m$-vectors such that 
$P'':\max\{ \bfc^\top \bfx \,|\, A \bfx \le \bfb'' \}$ and $P^*:\max\{ \bfc^\top \bfx \,|\, A \bfx \le \bfb^* \}$ are finite. 
Then, for each optimal solution $\bfx''$ of $P''$, there exists an optimal solution $\bfx^*$ of $P^*$ 
with $\| \bfx'' - \bfx^* \|_\infty \le n \Delta^* \| \bfb'' - \bfb^* \|_\infty$.	
\end{thm}
%\noindent
%Then, we have the following. 

\begin{cor}\label{lem:Reduction}
Let $\bfx''$ be an optimal solution of the rounded problem $(\ref{relaxP})$, and 
$J = \{ i \,|\, x''_i \ge n', \, i \in K \}$ as defined in Step 3 of the proposed algorithm. 
If the scaled problem $(\ref{scaleP})$ is feasible, the $i$-th coordinate of
the optimal solution of the scaled problem $(\ref{scaleP})$ is strictly positive 
for  $i \in J$.  Furthermore,  the same holds for the reduced problem $(\ref{smallP})$ 
and the original problem $(\ref{P})$ as the scaling factor $k$ is strictly positive. 
\begin{proof}
Define $\tilde{A} \in \bbR^{(2m'+n')\times n'}$, $\tilde{\bfb}''$, and $\tilde{\bfb}^* \in \bbR^{2m'+n'}$ as: 
\begin{eqnarray*}
\tilde{A} = 
\left[
\begin{array}{cc}
%I, \, (A_{L}')^{-1} A'_{\bar{L}}\\
%-I, \, -(A_{L}')^{-1} A'_{\bar{L}}\\
%-I
 E \\
- E \\
-I
\end{array}
\right]
, \,
\tilde{\bfb}^* = 
\left[
\begin{array}{c}
(A_{L}^{-1} \bfb')/k 	\\
-(A_{L}^{-1} \bfb')/k 	\\
\bfzero				
\end{array}
\right]
, \, 
\text{ and } \, 
\tilde{\bfb}'' = \lceil {\tilde{\bfb}^*} \rceil.  
\end{eqnarray*}
where $E=[I,\, (A_{L}')^{-1} A'_{\bar{L}}]$. 
With this notation, the rounded problem $(\ref{relaxP})$, respectively the scaled problem $(\ref{scaleP})$, 
can be restated as 
${P'':\max\{ -{\bfc'}^\top \bfx \,|\, \tilde{A} \bfx \le \tilde{\bfb}'' \}}$, 
respectively ${P^*:\max\{ -{\bfc'}^\top \bfx \,|\, \tilde{A} \bfx \le \tilde{\bfb}^* \}}$. 
Since $E$ is totally unimodular, $\tilde{A}$ is totally unimodular,  and thus  $\Delta^* = 1$ in Theorem~\ref{lem:Sc}. 
In addition, note that $\| \tilde{\bfb}'' - \tilde{\bfb}^* \|_\infty < 1$.
Recall that the scaled problem $(\ref{scaleP})$ and $P^*$ share the same unique optimal solution $\bfx^*$ 
as the optimal solution of the original problem $(\ref{P})$ is assumed to be unique. 
Therefore, since $\bfx''$ is an optimal solution of $P''$, we observe that 
$\| \bfx'' - \bfx^* \|_\infty < n'$ by Theorem~\ref{lem:Sc} and thus,  $x^*_i > 0$ for $i \in J$. 
\end{proof}
\end{cor}

Finally, we show the strong polynomiality of the proposed algorithm using Kitahara and Mizuno~\cite{KiMi12, KiMi13} results showing that 
the number of different basic feasible solutions generated by the primal simplex method 
with the most negative pivoting rule -- Dantzig's rule -- or the best improvement pivoting rule is bounded by:
\begin{equation*}
n \lceil m \frac{\gamma}{\delta} \log(m \frac{\gamma}{\delta})  \rceil
\end{equation*}
where $m$ is the number of constraints, $n$ is the number of variables, 
and $\gamma$ and $\delta$ are, respectively,  the minimum and the maximum values 
of all the positive elements of the primal basic feasible solutions.
Thus, we need to estimate the values $\gamma$ and $\delta$ for the introduced auxiliary problems.

\par
Since the coefficient matrices used in the proposed  algorithm are totally unimodular and the right hand side vector of 
the rounded problem $(\ref{relaxP})$ %and $(\ref{fsblP})$ are integers, 
is integer, we have $\delta=1$.
For $\gamma$, we use Lemma~\ref{lem:gamma}. 
\begin{lem}\label{lem:gamma}
For the auxiliary problem $(\ref{relaxP})$, %  and  $(\ref{fsblP})$, 
we have $\gamma \le \gamma^* = m( {m}{n}(m + n^2)+1)$. 
\begin{proof}
Note that the right-hand side vector for  $(\ref{relaxP})$ %and $(\ref{fsblP})$ are both 
is $\lceil (A_{L}')^{-1} \bfb' / k \rceil$. 
By the total unimodularity, we observe that 
\begin{equation*}
\| \lceil (A_{L}')^{-1} \bfb'/k \rceil \|_\infty \le \| (A_{L}')^{-1} \bfb'/k \|_\infty + 1 \le m' \| \bfb' \|_\infty /k + 1.
\end{equation*}
The numerator $\| {A'}^\top ({A'}{A'}^\top)^{-1} \bfb' \|_{2}$ of $k$ is bounded below by $\| \bfb' \|_\infty / n'$
implying $\| \lceil (A_{L}')^{-1} \bfb'/k \rceil \|_\infty \le {m'}{n'}(m' + (n')^2)+1$. 
Thus, by Cramer's rule and the total unimodularity of the coefficient matrix of $(\ref{relaxP})$, % and $(\ref{fsblP})$, 
the $l_\infty$-norm of a basic solution of $(\ref{relaxP})$ %or $(\ref{fsblP})$ 
is bounded above by $m'( {m'}{n'}(m' + (n')^2)+1)$. 
\end{proof}
\end{lem}
\noindent
The two-phase simplex algorithm is called at most $m$ times. Thus, the number of auxiliary problems solved 
by the proposed algorithm is bounded above by $2m$ as each call corresponds to 2 auxiliary problems : one for each phase.
%when we count two problems for .  
%and the auxiliary problem $(\ref{fsblP})$ is solved at most $m$ times for each $(\ref{relaxP})$. 
Therefore, if all the auxiliary problems are non-degenerate,  the total number of basic solutions generated by the 
algorithm is bounded above by  $2m [ n \lceil m \gamma^* \log(m \gamma^*)  \rceil ]$; that is by
\begin{equation*}
 2mn \lceil (m^4 n + m^3 n^3 + m^2) \log (m^4 n + m^3 n^3 + m^2) \rceil 
\end{equation*}
which completes the proof of Theorem~\ref{thm1}. Alternatively, since $m \le n$, this bound can be restated as $O(m^4 n^4 \log n)$. 
While assuming the non-degeneracy of the auxiliary problems is needed to use Kitahara-Mizuno's bound, 
the number of degenerate updates of bases at a single basic solution is typically not too large in practice. 

%\par
%It should be noted that Mizuno~\cite{Mi14} also estimates the bound on the total number of basic solutions generated by his algorithm. 
%For totally unimodular LP problems, his bound can be simplified to $O(m^3n^4 \log( m n^3))$. 
%Therefore, in a theoretical manner, one could say that our algorithm is inferior to that of Mizuno~\cite{Mi14}. 
%However, it would be interesting to compare their performances through numerical experiments. 
%For instance, one could expect that our algorithm outperforms for problems with relatively small $m$, 
%because it stops after identifying at most $m$ coordinates strictly positive at optimality, 
%while his algorithm needs to identify $(n-m)$ coordinates equal zero for the worst case. 

%\par
%Assuming the nondegeneracy is not strong in practice, 
%because it is seldom to have many degenerate updates of bases  
%at a single basic solution. 

\vspace{5mm}
\noindent{\bf Acknowledgment}\\
Research supported in part by Grant-in-Aid for Science Research (A) 
26242027 of Japan Society for the Promotion of Science (JSPS), 
Grant-in-Aid for JSPS Fellows,
the Natural Sciences and Engineering Research Council of Canada, the Digiteo Chair
C\&O program, and by the Canada Research Chairs program.
Part of this research was done while the authors were at the LRI, Universit\'e Paris-Sud, Orsay, France, within the Digiteo invited researchers program.

%%%

%%%

\end{document}